\documentclass{amsart}

\newtheorem{theorem}{Theorem}[section]
\newtheorem{lemma}[theorem]{Lemma}
\newtheorem{proposition}[theorem]{Proposition}

\newtheorem*{claim}{Claim}

\theoremstyle{definition}

\newtheorem{example}[theorem]{Example}
\newtheorem{remark}[theorem]{Remark}
\newtheorem*{acknowledgement}{Acknowledgement}

\theoremstyle{remark}

\numberwithin{equation}{section}

\usepackage{amscd, amssymb}
\input{diagrams}

\newcommand{\NN}{\mathbb{N}}
\newcommand{\ZZ}{\mathbb{Z}}
\newcommand{\QQ}{\mathbb{Q}}
\newcommand{\RR}{\mathbb{R}}

\newcommand{\PP}{\mathbb{P}}

\newcommand  {\shE}     {\mathcal{E}}
\newcommand  {\shF}     {\mathcal{F}}

\newcommand  {\shM}     {\mathcal{M}}

\newcommand  {\shL}     {\mathcal{L}}

\newcommand  {\foS}     {\mathfrak{S}}

\newcommand  {\aff}     {{\text{aff}}}

\newcommand  {\Aut}     {\operatorname{Aut}}

\renewcommand{\cong}    {\equiv}

\newcommand  {\co}       {\colon}

\newcommand  {\Div}     {\operatorname{Div}}

\newcommand  {\Fr}      {\operatorname{Fr}}

\renewcommand  {\k}     {\kappa}

\newcommand  {\lra}     {\longrightarrow}

\newcommand  {\N}       {\operatorname{N}}

\newcommand  {\NEP}     {\operatorname{\overline{NE}}}

\newcommand  {\NS}      {\operatorname{NS}}

\renewcommand{\O}       {\mathcal{O}}

\renewcommand{\P}       {\overline{P}}

\newcommand  {\Pic}     {\operatorname{Pic}}

\newcommand  {\Proj}    {\operatorname{Proj}}

\newcommand  {\ra}      {\rightarrow}

\newcommand  {\rank}    {\operatorname{rank}}
\newcommand  {\red}     {{\operatorname{red}}}
\newcommand  {\Reg}     {\operatorname{Reg}}

\newcommand  {\SBs}     {\operatorname{SBs}}

\newcommand  {\SL}      {\operatorname{SL}}
\newcommand  {\Spec}    {\operatorname{Spec}}
\newcommand  {\Supp}    {\operatorname{Supp}}

\begin{document}

\title[Normal del Pezzo surfaces]{Normal del Pezzo surfaces containing a nonrational singularity}

\author[Stefan Schroeer]{Stefan Schr\"oer}

\address{Mathematische Fakult\"at, Ruhr-Universit\"at, 
         44780 Bochum, Germany}

\curraddr{M.I.T. Mathematical Department, 
          77 Massachusetts Avenue, Cambridge, MA 02139-4307, USA}

\email{s.schroeer@ruhr-uni-bochum.de}


\dedicatory{Revised version, 28 Oktober 2000}

\begin{abstract}
Working over perfect ground fields of arbitrary characteristic, I
classify minimal normal del Pezzo surfaces containing a 
nonrational singularity.
As an application, I determine the structure of 2-dimensional 
anticanonical models for proper  normal algebraic surfaces.
The anticanonical ring may be non-finitely generated.
However, the anticanonical model is either a proper surface, 
or a proper surface minus a point.
\end{abstract}

\maketitle

\section*{Introduction}

Proper smooth algebraic surfaces with ample anticanonical divisor are called 
del Pezzo surfaces. Prominent examples are 
$\PP^2$ and 
$\PP^1\times\PP^1$. We call a  proper \emph{normal} algebraic surface
$X$  del Pezzo if  
$(-K_X)^2>0$ and 
$(-K_X)\cdot C> 0$ holds for all curves 
$C\subset X$. Here we use
Mumford's rational intersection numbers.
Note that $X$ might be non-$\QQ$-Gorenstein, such that the antipluricanonical 
sections do not define a closed embedding into any projective space.

Several authors studied normal del Pezzo surfaces. For example, 
Hidaka and Watanabe \cite{Hidaka;
Watanabe 1981} considered the Gorenstein case. Sakai \cite{Sakai 1984} 
studied complex rational
surfaces. B\u{a}descu \cite{Badescu 1983} analyzed complex ruled surfaces. 
Fujisawa  \cite{Fujisawa 1995} and Cheltsov \cite{Cheltsov 1997} classified  
complex normal del Pezzo
surfaces.

Here we shall study normal del Pezzo surfaces
over perfect ground fields of arbitrary characteristic. Contractions of 
sections on suitable
$\PP^1$-bundles   give plenty of  examples.
The main result is that, starting from such contracted bundles, we can reach 
any 
normal del Pezzo surface with base number 
$\rho(X)=1$ containing a nonrational singularity via a sequence of   
 \emph{generalized elementary transformations}.     Such elementary 
transformations are closely related
to  the  monoid 
$\SL_2(\NN)$  and continued fractions.

We can view such del Pezzo surfaces
as ``minimal models''  in the category of normal surfaces. 
As an application of our results, we determine the structure
of anticanonical models 
$P(-K_X)=\Proj (\bigoplus_{n\geq 0} H^0(X,-nK_X))$ for normal surfaces.
It turns out that each  2-dimensional anticanonical model  is either a 
proper normal $\QQ$-Gorenstein surface, or there is a
canonical compactification 
$P(-K_X)\subset \P$ by adding a single non-$\QQ$-Gorenstein point at infinity.
In the latter case, the anticanonical ring
$R(-K_X)=\bigoplus_{n\geq 0} H^0(X,-nK_X)$
is not finitely generated. This happens, for example, if
the stable base locus of $-K_X$ contains isolated non-$\QQ$-Gorenstein
singularities.

\begin{acknowledgement}
I wish to thank the referee for his remarks and suggestions, which
helped to clarify the paper.
\end{acknowledgement}

\section{Surfaces with ineffective canonical class}

We fix some notation.
Throughout the paper, we    work over a perfect ground field 
$k$ of arbitrary  characteristic 
$p\geq 0$. Suppose
$X$ is a proper normal algebraic surface with $k=\Gamma(X,\O_{X})$. Let
$p_g(X)=h^2(\O_{X})$ be  the
\emph{geometric genus} and 
$q(X)=\dim(\Pic^0_X)$   the \emph{irregularity}. We are mainly 
interested in surfaces with $p_g(X)=0$. For such surfaces, the Picard
scheme is smooth \cite{Mumford 1966} p.\ 198, so the irregularity
is also $q(X)=h^1(\O_X)$.

Two Weil divisors 
$A, B\in Z^1(X)$ are   \emph{numerically equivalent}, written 
$A\equiv B$, if 
$A\cdot C = B\cdot C$ for all curves 
$C\subset X$. Here we 
use Mumford's  rational intersection numbers \cite{Mumford 1961}.
Numerical equivalence yields  two quotients 
$Z^1(X)\ra N(X)$ and
$\Pic(X)\ra N^1(X)$, together with an inclusion
$N^1(X)\subset N(X)$. The rank of 
$N^1(X)$ is called the \emph{base number} 
$\rho(X)$.
Note that the base number might be zero \cite{Schroeer 1999a}. However, 
the following holds.

\begin{proposition}
\label{Q-Cartier}
Suppose 
$p_g(X)=0$. Then each Weil divisor on 
$X$ is numerically equivalent to a 
$\QQ$-Cartier divisor. Moreover,
$X$ is projective.
\end{proposition}

\proof
See \cite{Schroeer 1999b} Corollary 4.4.
\qed

\medskip
Let 
$f\co  X\ra Y$ be a proper morphism with
$\O_Y= f_*( \O_X)$ onto another proper normal algebraic scheme $Y$. 
If $Y$ is a curve, we  say that 
$f\co X\ra Y$ is a \emph{fibration}. The arithmetic genus 
$p_a(X_\eta)=h^1(\O_{X_\eta})$ of the generic fiber is called the 
\emph{genus of the fibration}.
If 
$Y$ is a surface, we say that
$f\co X\ra Y$ is a (birational) \emph{contraction}. A curve 
$R\subset X$ is called \emph{contractible} if there is a contraction 
$f\co X\ra Y$ such that 
$R$ is the exceptional curve.

\begin{proposition}
\label{contraction on fibration}
Assume that
$p\co X\ra B$ is a fibration of genus zero. Then 
$R^1p_*(\O_X)=0$. Moreover, for  all
$b\in B$,  each   curve 
$R\varsubsetneq p^{-1}(b)$ is contractible.
\end{proposition}

\proof
Passing to a resolution of singularities and a relatively minimal model, we  
easily deduce
$R^1p_*(\O_X)=0$. Now the contractibility of $R$ follows   from Artin's 
contraction
criterion (\cite{Artin 1962} Thm.~2.3).
\qed

\medskip
The \emph{pseudoeffective cone} $\NEP(X)\subset N(X)_\RR$ is the closed 
convex cone 
generated by all curves 
$C\subset X$.  The pseudoeffective cone
has full dimension and is generated by its extremal rays. 
According to Koll\'ar \cite{Kollar 1995}, Lemma 4.12, each extremal ray 
$\RR_+E\subset \NEP(X)$ with 
$E^2<0$ is generated by an integral curve 
$R\subset X$. 

Suppose the pseudoeffective cone is simplicial, in other words,  the 
convex hull
of 
$n=\rank N(X)$ linearly independent extremal rays 
$\RR_+E_i\subset \NEP(X)$. This holds for example if $\rank N(X)=2$. 
Decompose the canonical class
$$
K_X\equiv \lambda_1E_1+\lambda_2E_2 +\ldots +\lambda_nE_n,
$$
for certain  unique coefficients 
$\lambda_i\in\RR$.

\begin{proposition}
\label{contractible curve}
With the preceding assumptions, suppose 
$E_i^2<0$ and 
$\lambda_j<0$ for some pair $j\neq i$. Then the curve 
$R_i\subset Y$ generating  
$\RR_+E_i$ is contractible. Moreover, if  
$f\co X\ra Y$ is the corresponding contraction, then   
$K_Y$ is not pseudoeffective.
\end{proposition}

\proof
The conditions ensure that 
$K_X+nR_i$ is not pseudoeffective for all integers 
$n\geq 0$. According to \cite{Schroeer 1999b}, Corollary 5.2, the contraction 
$f\co X\ra Y$ of 
$R$ exists. Seeking a contradiction, assume that  
$K_Y$ is pseudoeffective. Then  
$f^*(K_Y)$ is also pseudoeffective. Write 
$K_X=f^*(K_Y)+\lambda R_i$ for some coefficient 
$\lambda\in\QQ$.  Then $K_X+nR_i$  is   pseudoeffective for some 
$n\gg 0$, a contradiction. Thus 
$K_Y$ is not pseudoeffective.
\qed

\medskip
A proper normal algebraic surface with 
$(-K_X)^2>0$ and 
$(-K_X)\cdot C> 0$  for all curves 
$C\subset X$ is called a \emph{del Pezzo surface}.
This holds, for example, if  $-K_X$ is an ample 
$\QQ$-Cartier divisor. 

\begin{proposition}
\label{contraction of del Pezzo}
The contraction of a del Pezzo surface is a del Pezzo surface.
\end{proposition}

\proof
Let 
$f\co X\ra Y$ be such a contraction. According to the projection formula,
$K_Y\cdot C<0$ for all curves 
$C\subset Y$. Writing 
$K_X=f^*(K_Y)+ K_{X/Y}$ with 
$K_{X/Y}$ supported on the exceptional curve, we obtain 
$$
0<K_X^2 = K_Y^2 + 2f^*(K_Y)\cdot K_{X/Y} + K_{X/Y}^2 = 
K_Y^2  + K_{X/Y}^2\leq K_Y^2.
$$
Therefore 
$Y$ is del Pezzo.
\qed

\medskip
Recall that a curve 
$R\subset X$ is called \emph{negative definite} if  the intersection matrix 
$(R_i\cdot R_j)$ is negative definite, where
$R_i\subset R$ are the irreducible components. Contractible curves 
are negative definite.
For   surfaces with many antipluricanonical sections, the converse holds 
as well: 

\begin{proposition}
\label{contractions on del Pezzo}
Suppose there is an integer 
$m>0$ with 
$h^0(-mK_X)>1$. Then each negative definite curve 
$R\subset X$  is contractible.
\end{proposition}

\proof
By induction on the number of irreducible components of 
$R$, it suffices to treat the case that 
$R$ is integral.
Suppose
$D=K_X+nR$ is  effective for some integer 
$n\geq 0$. Then  
$-mK_X= mnR-mD$, so we have an inclusion 
$H^0(-mK_X)\subset H^0(mnR)$. Using 
$R^2<0$, you easily see that 
$h^0(mnR)=1$, contradiction. 
Consequently,  
$K_X+nR$ is not effective for all 
$n\geq 0$.
Now \cite{Schroeer 1999b}, Corollary 5.2 ensures that 
$R\subset X$ is contractible.
\qed

\medskip
The next result reduces all surfaces with 
ineffective canonical class  to a  special situation:

\begin{theorem}
\label{reduction}
Suppose 
$X$ is a proper normal algebraic surface with 
$-K_X$ not pseudoeffective. Then there is a contraction 
$h\co X\ra Y$ such that one of the following holds:
\renewcommand{\labelenumi}{(\roman{enumi})}
\begin{enumerate}
\item  The contraction
$Y$ is a del Pezzo surface with base number 
$\rho(Y)=1$.
\item There is a genus zero  fibration 
$p\co Y\ra B$ with irreducible fibers, and the base number is
$\rho(Y)=2$.
\end{enumerate}
\end{theorem}

\proof
For 
$\rho(X)=1$ there is nothing to prove. Suppose 
$\rho(X)\geq 2$. Choose a pseudoample class 
$A\in N(X)_\RR$ with 
$K_X\cdot A<0$, and let 
$\RR_+E\subset \NEP(X)_\RR\cap A^\perp$ be an extremal ray. 

First, suppose 
$E^2<0$. Then $\RR_+E$ is generated by an integral curve 
$E\subset X$. Since 
$(K_X+nE)\cdot A<0$ holds for all 
$n\geq 0$, the class 
$K_X+nE$ is not pseudoeffective, so the contraction 
$X\ra X'$ of 
$B$ exist by \cite{Schroeer 1999b}, Corollary 5.2.
Since $A$ comes from an ample class on 
$X'$, the class
$K_{X'}$ is not pseudoeffective, and we can proceed by induction on the 
base number
$\rho(X)$.

Second, suppose 
$E^2=0$. Then the existence of a genus zero fibration 
$p\co X\ra B$  easily follows from Mori 
\cite{Mori 1982}, Theorem 2.3.
Since 
$\RR_+ E$ is extremal, each fiber 
$X_b$, 
$b\in B$ is irreducible.
Using  
$\Pic^0(Y_\eta)=0$, we directly deduce  
$\rho(Y)=2$.
\qed

\begin{remark}
If 
$p\co Y\ra B$ is a genus zero fibration  with irreducible fibers, 
it may or may not be possible to contract
$Y$, and the result may or may not be del Pezzo. We take up this 
issues in Sections 
\ref{contractions of ruled} and \ref{Classification of del Pezzo surfaces}.
\end{remark}

\section{Surfaces containing  a nonrational singularity}

The aim of this section 
is to determine the structure of proper normal algebraic  surfaces 
$Z$ with ineffective pluricanonical  classes containing  a 
nonrational singularity. 
The following   observation will be useful:

\begin{lemma}
\label{torsion-free}
Suppose 
$Z$ is a proper normal algebraic surface with irregularity
$q(Z)=0$. If no multiple 
$nK_Z$ with 
$n> 0$ is   effective, then 
$\Pic( Z)$ is a   free 
$\ZZ$-module of finite rank.
\end{lemma}

\proof
Assume there is an invertible 
$\O_Z$-module 
$\shL\in\Pic(X)$ of order 
$n>0$. The Riemann--Roch theorem tells us 
$\chi(\shL)=\chi(\O_Z) >0$. On the other hand, 
$$
\chi(\shL)\leq h^0(Z,\shL) + h^0(Z,\shL^{\otimes -1}\otimes\omega_Z).
$$
Since 
$Z$ is integral and 
$\shL$ is nontrivial, the first summand is zero. Because the invertible sheaf
$\shL^{\otimes -n}(nK_Z) =\O_Z(nK_Z)$  is not effective, the 
second summand vanishes as well.
Thus 
$\chi(\shL)\leq 0$, a contradiction. Hence 
$\Pic(Z)$ is torsion free. 

The group scheme
$\Pic^0_Z$ vanishes  because its Lie algebra 
$H^1(Z,\O_{Z})=0$ does, so the map 
$\Pic(Z)\ra \NS(Z)$ onto the N\'eron--Severi group is bijective. 
Since the N\'eron--Severi
group is finitely generated,
$\Pic(Z)$ is free and finitely generated.
\qed

\medskip
Suppose 
$Z$ is a proper normal algebraic  surface. Let 
$f\co X\ra Z$ be a resolution of singularities, and
$R\subset X$ the exceptional divisor.
Choose a maximal subcurve 
$R'\subset R$ whose formal completion 
$X_{/R'}$ satisfies 
$H^1(\O_{X_{/R'}})=0$. According to Artin \cite{Artin 1962}, Theorem 2.3, 
the contraction 
$h:X\ra Y$ of 
$R'\subset X$ exists. We call the corresponding partial resolution 
$g:Y\ra Z$ a \emph{minimal resolution of nonrational singularities}.

\begin{theorem}
\label{structure}
Suppose  
$Z$ is a proper normal algebraic  surface 
 containing a nonrational singularity, and that no multiple 
$nK_Z$ with
$n>0$ is   effective. Let 
$\co Y\ra Z$ be a minimal resolution of nonrational singularities and 
$S\subset Y$  the exceptional curve.
Then the following hold:
\renewcommand{\labelenumi}{(\roman{enumi})}
\begin{enumerate}
\item 
There is a genus zero  fibration 
$p\co Y\ra B$  over a curve 
$B$ of genus 
$g>0$.
\item
The curve 
$S\subset Y$ is a section over 
$B$, and each singularity on 
$Z$ is rational except for
$z=g(S)$. 
\item 
The irregularity is
$q(Z)=0$, and 
$\Pic(Z)$ is a free 
$\ZZ$-module of finite rank.
\end{enumerate}
\end{theorem}

\proof
Let $f\co X\ra Z$ be the minimal resolution of singularities.
Clearly, $nK_X\neq 0$ for all $n>0$.
Suppose there is a nonzero section $s\in H^0(X,nK_X)$.
The corresponding curve is supported on $R$, and this implies
$K_X\cdot R_i<0$ for some irreducible component $R_i\subset R$, 
contradicting the minimality
of the resolution. Therefore $X$ has Kodaira dimension  $\kappa(X)=-\infty$.
The Enriques--Kodaira classification ensures the existence of  a fibration 
$\overline{X}\ra \overline{B}$ of genus zero    over the algebraic closure 
$k\subset \bar{k}$ 
(see Mumford \cite{Mumford 1969} Sect.\ 1). Since 
$p_g(Z)=0$, the sequence 
$$
0 \lra H^1(Z,\O_Z)  \lra H^1(X,\O_X)  \lra H^1(X_{/R},\O_{X_{/R}})  \lra 0
$$
is exact, where 
$R\subset X$ is the exceptional divisor and 
$X_{/R}$ is the corresponding formal completion. As 
$H^1(\O_{\overline{B}})\ra H^1( \O_{\overline{X}})$ is bijective, the  curve 
$\overline{B}$ has genus
$g>0$. Thus each horizontal curve on
$ \overline{X}$ has positive  arithmetic genus. Consequently, the Galois group 
$\Aut(\bar{k}/k)$ permutes the vertical curves, 
so the fibration descends to a fibration 
$X\ra B$. The exceptional curve 
$R'\subset X$ for the contraction 
$h\co X\ra Y$
must be vertical, hence there is an induced fibration 
$p\co Y\ra B$. This proves  (i).

By Proposition \ref{contraction on fibration}, the exceptional divisor 
$S\subset Y$ is horizontal. For each integer 
$n\geq 0$, the kernel of the restriction homomorphism 
$\Pic^0_{nS}\ra\Pic^0_S$ is a smooth unipotent group scheme. 
On the other hand, 
$\Pic^0_B$ is an Abelian variety. The restriction map 
$\Pic^0_B\ra\Pic^0_{nS}$ is an epimorphism, since the induced map 
$H^1(\O_B)\ra H^1(\O_{nS})$ on Lie algebras is surjective. We deduce 
$\Pic^0_{nS}=\Pic^0_S$ and 
$H^1(\O_{Y_{/S}})=H^1(\O_S)$, where 
$Y_{/S}$ is the formal completion of $S\subset Y$.
This gives an  exact sequence 
$$
0 \lra \Pic^\tau_Z  \lra \Pic^0_B  \lra \Pic^0_S  \lra 0
$$
of group schemes. Here 
$\Pic^\tau_Z $ is the group scheme of numerically trivial invertible sheaves. 
Since 
$ \Pic^0_B $ is an  abelian variety, the kernel  
$\Pic^\tau_Z $ is proper. According to  \cite{SGA 6} Sect.~XII, Cor.~1.5, 
it  must be affine, because 
$S\ra B$ is proper.
 Since 
$H^2(\O_{Z})=0$, the group scheme 
$\Pic^\tau_Z$ is smooth (see \cite{Mumford 1966}  p.~198).
Being smooth, proper, and affine,
$\Pic^\tau_Z$ is \'etale. This implies 
$q(Z)=0$. According to Lemma \ref{torsion-free},
$\Pic(Z)$ is free and finitely generated, so $\Pic^\tau_Z=0$. 
This proves  (iii).

It remains to verify  (ii). Since 
$S\subset Y$ is horizontal, each integral component  of $S$ has genus 
$g>0$. Hence 
$S$ is  geometrically integral, and  in particular 
$k\subset\Gamma(\O_S)$ is bijective. We also see that   
$z=g(S)$ is the only nonrational singularity on the surface
$Z$. 
Assume that
$S\subset Y$ is not a section.
Using that 
$\Pic_B^0\ra\Pic^0_S$ is an isomorphism, we  derive a contradiction as follows.

First, assume that the normalization  
$\tilde{S}\ra S$ is not an isomorphism. Since 
$S$ is irreducible, the kernel of 
$\Pic_S^0\ra\Pic_{\widetilde{S}}^0$ is a  
nonfinite  unipotent  group scheme, contradicting
that 
$\Pic_S^0=\Pic_B^0$ is abelian. Consequently, 
$S$ is normal.
Next, decompose 
$S\ra B$ into a purely inseparable morphism 
$S\ra A$ and a separable morphism 
$A\ra B$. Using 
$\Pic_B^0=\Pic^0_S$, we deduce that 
$\Pic_B^0\ra\Pic^0_A$ and 
$\Pic_A^0\ra\Pic^0_S$ are isomorphisms. Set 
$n=\deg(A/B)$. The Hurwitz formula gives 
$$
2g-2 = n(2g-2) +\deg(K_{A/B}) \geq n(2g-2).
$$
If $n>1$, then
$g=1$, and 
$A\ra B$ is an isogeny of degree 
$n$ between elliptic curves. So the dual isogeny  
$\Pic_B^0\ra\Pic^0_A$ has a kernel of length 
$n$,  hence 
$n=1$. Consequently, the projection
$S\ra B$ is purely inseparable.
Finally, set 
$p^m=\deg(S/B)$. The two curves $A$ and $B$
 are isomorphic without 
$k$-structures, and the projection 
$S\ra B$ is just the iterated linear Frobenius map
$\Fr^m\co B\ra B$. The induced morphism 
$\Fr^*\co \Pic_B^0\ra\Pic_B^0$ is   multiplication by 
$p$, which  has a kernel of length 
$p^{2g}$ (see Mumford \cite{Mumford 1970} p.64). Hence 
$m=0$, and 
$S\subset Y$ is a section.
\qed

\begin{remark}
In the preceding result, 
the generic fiber of the genus zero fibration 
$Y\ra B$ is a form of the projective line.  However, since 
$S\subset Y$ is a section, 
$Y_\eta$ contains a rational point and is therefore  a   projective line
over the function field 
$\kappa(B)$.
\end{remark}

\begin{remark}
The   arguments in this section do not work over nonperfect ground fields. 
Curves with ``genus change" cause problems.
Compare Tate  \cite{Tate 1952}.
\end{remark}

\section{Contractions of $\PP^1$-bundles}\label{contractions of ruled}

The simplest   del Pezzo surfaces containing a nonrational singularity
are contractions of 
$\PP^1$-bundles. The task
now is to   discuss such contractions. 
Fix a smooth curve 
$B$ of genus 
$g>0$ and a 
$\PP^1$-bundle 
$q\co Y\ra B$. Then the pseudoeffective cone
$\NEP(Y)\subset \N(Y)_\RR$ is 2-dimensional and generated by two extremal rays.
The first extremal ray is the  \emph{fiber class} 
$F\in N(Y)_\QQ$, defined by the equations 
$Y_b = \dim_k\kappa(b) \cdot F$ for each 
$b\in B$. The second extremal ray may or may not be generated by a  curve.

\begin{lemma}
\label{canonical class}
Let 
$R\subset Y$ be a horizontal curve, say of degree 
$n$ over 
$B$. Then 
$n$ divides $R^2$, and 
$nK_{Y/B}\cong -2R + R^2/n\cdot F$.
\end{lemma}

\proof
Write 
$nK_{Y/B}\cong-2R+\lambda F$ for some   
$\lambda\in\ZZ$. Using $K_{Y/B}^2=0$, we obtain
$0= (nK_{Y/B})^2 = 4R^2 -4\lambda n$,
hence the result.
\qed

\medskip
Set 
$e=-\inf \left\{ S^2\mid S\subset Y \text{ a section} \right\}$. 
Note that there is a section 
$S\subset Y$ with 
$S^2<0$ if and only if 
$e>0$.

\begin{proposition}
\label{del Pezzo contraction}
The surface $Y$ contracts to a normal del Pezzo surface if and only if 
$2g-2-e<0$.
\end{proposition}

\proof
The condition is sufficient: Let 
$S\subset Y$ be the section with 
$S^2=-e<0$. Since 
$K_Y\equiv-2S+(2g-2-e)F$, the contraction exists by Proposition 
\ref{contractible curve}. Since 
$K_Z=(2g-2-e)\cdot g_*(F)$, the resulting surface is del Pezzo.
The condition is necessary: By Theorem \ref{structure}, the exceptional curve 
$S\subset Y$ is a section. Since 
$2g-2-e$ is the multiplicity of 
$K_Z$, the condition follows.
\qed

\begin{proposition}
\label{no separable curve}
Suppose 
$R\subset Y$ is an integral curve with 
$R^2<0$. Then the projection 
$q\co R\ra B$ is purely inseparable.
\end{proposition}

\proof
The separable closure of the field extension $\kappa(B)\subset\kappa(R)$
defines a proper normal algebraic curve $A$, together
with a purely inseparable morphism $R\ra A$ and a separable morphism
$A\ra B$. Set $n=\deg(A/B)$. On the induced $\PP^1$-bundle $Y_A$,
the preimage of $R$ splits into $n$ irreducible  components
$R_A=R_1+\ldots+R_n$, which are permuted by the Galois action.
This implies $R_1^2=\ldots =R_n^2$, hence all irreducible components
have $R_i^2<0$. Since the pseudoeffective cone $\NEP(Y_A)$
has at most one extremal ray with negative self intersection,
we conclude $n=1$. Thus $A=B$, and $R\ra B$ is purely inseparable.
\qed

\medskip
In characteristic 
zero, we conclude that each integral curve 
$R\subset Y$ with 
$R^2<0$ is necessarily a section. 
Suppose there is such a section. Setting 
$\shE=p_*(\O_Y(R))$ and 
$\shL=p_*(\O_S(R))$, we obtain an extension 
$$
0 \lra \O_B  \lra \shE  \lra \shL  \lra 0
$$
with 
$X=\PP(\shE)$ and 
$R=\PP(\shL)$.

\begin{proposition}
\label{contractible for p=0}
Suppose the characteristic is
$p=0$, and let 
$R\subset Y$  be as above. Then 
$R$ is contractible if and only if the extension 
$\shE$ splits.
\end{proposition}

\proof
If there is a splitting 
$\shE\ra\O_B$, then 
$A=\PP(\O_B)$ defines another section disjoint from 
$S$ with $A^2>0$. By the 
Fujita--Zariski Theorem (\cite{Fujita 1983} Thm.~1.10),  
$\O_Y(A)$ is pseudoample, hence a suitable multiple is base point free and 
defines the desired contraction.

Conversely, assume that the contraction 
$g\co Y\ra Z$  exists. A direct calculation gives
$$
\Pic(2R)=\Pic(R)\oplus H^1(S,\O_R(-R))
$$
Since 
$Z$ contains only one singularity, it is projective.
This is because the complement of an affine neighborhood
of the singularity supports an ample Cartier divisor,
as explained in \cite{Hartshorne 1970} p.\ 69.
Let 
$\shL$ be the preimage of an ample invertible 
$\O_Z$-module. We calculate the restriction 
$\shL|_{2R}$ in two ways. Since 
$\shL$ comes from 
$Z$, the restriction is trivial. On the other hand, $\shL|_{2R}$ defines a 
class
$  (0,\alpha)\in\Pic(2R)$ according to the
decomposition of
$\Pic(2R)$. A straightforward cocycle computation left to the reader reveals
that 
$\alpha\in H^1(R,\O_R(-R))$   coincides with a multiple of the 
extension class of
$\shE$ in 
$H^1(B,\shL^\vee)$. Since the characteristic is 
zero, this extension class   is zero.
\qed

\medskip
The situation in positive characteristic is different:

\begin{theorem}
\label{contractible for p>0}
Suppose the characteristic is
$p>0$. Let 
$\shE$ be a rank two vector bundle on $B$ with 
$X=\PP(\shE)$. Then the following are equivalent:
\renewcommand{\labelenumi}{(\roman{enumi})}
\begin{enumerate}
\item
For some integer 
$m\geq 0$, there is a splitting 
$(\Fr^m)^*\shE=\shL_1\oplus\shL_2$ with invertible summands satisfying 
$\deg(\shL_1)\neq\deg(\shL_2)$.
\item 
There is a contractible curve 
$R\subset Y$.
\item
There is a curve 
$R\subset Y$ with 
negative self-intersection.
\end{enumerate}
\end{theorem}

\proof
First, note that 
$\PP(\shF)\simeq\PP(\shE)$ if and only if there is an invertible sheaf 
$\shL$ with 
$\shF\simeq \shE\otimes\shL$. Hence the condition in (i) does not depend 
on the  choice of
$\shE$.

Suppose (i) holds, say with 
$\deg(\shL_1)<\deg(\shL_2)$. Set 
$Y_m=\PP((\Fr^m)^*\shE)$. Then there is a cartesian diagram 
$$
\begin{CD}
Y_m     @>>> Y    \\
@VVV     @VVV\\
B    @>>\Fr^m>   B
\end{CD}
$$
whose horizontal maps are bijective. The section 
$R_m=\PP(\shL_1)$ of $Y_m$ satisfies 
$(R_m)^2=\deg(\shL_1)-\deg(\shL_2)<0$. Since $\PP(\shL_2)$ is another 
disjoint section, we conclude as
in the proof of Proposition
\ref{contractible for p=0} that 
$R_m\subset Y_m$ is contractible. Similarly, the  image 
$R\subset Y$ of 
$R_m$ is contractible as well. Hence (ii) holds.

The implication (ii) $\Rightarrow$ (iii) is well known.
Finally, suppose an integral curve 
$R\subset Y$ has 
$R^2<0$. By Proposition \ref{no separable curve}, the projection 
$R\ra B$ is purely inseparable, say of degree 
$p^m$. Hence its  normalization $\tilde{R}$ is isomorphic to
$B$, and the induced morphism 
$\tilde{R}\ra B$ it the iterated linear Frobenius
$\Fr^m\co B\ra B$. Making base change along this morphism, we can assume that 
$R\subset Y$ is a section. The section defines an extension 
$$
0 \lra \O_B  \lra \shE   \lra \shL  \lra 0
$$
with 
$Y=\PP(\shE)$ and 
$\shL=p_*(\O_R(R))$. 
The  extension class for $\shE$ lies in 
$H^1(B,\shL^{-1})$. Since
$\deg(\shL^{-1})= - R^2>0$, the sheaf
$\shL^{-1}$ is ample, so  $H^1(B,\shL^{-p^m})=0$ for all
$m$ sufficiently large. Since $\shL^{-p^m} =(\Fr^m)^*\shL^{-1}$, 
the induced extension
$$
0 \lra \O_B  \lra (\Fr^m)^*\shE   \lra (\Fr^m)^*\shL  \lra 0
$$
splits. Thus condition (i) holds.
\qed

\begin{remark}
The existence of nonsplit rank two vector bundles that split
on taking a purely inseparable cover follows from Tango's work.
In \cite{Tango 1972}, Theorem 15, he showed that on 
smooth prper curves $B$ satisfying certain conditions related
to the Hasse--Witt matrix, there  is an ample invertible $\O_X$-module
$\shM$ so that the Frobenius map 
$\Fr\co H^1(B,\shM^\vee)\ra H^1(B,\Fr^*\shM^\vee)$ is not injective. Then
the corresponding extension $0\ra\O_B\ra\shE\ra\shM\ra 0$ 
gives the desired rank two vector
bundle $\shE$. For explicit examples, see \cite{Tango 1972}, Section 5.
\end{remark}

\section{Elementary transformations and continued fractions}
\label{Elementary transformations}

In this section, we calculate the behavior of   multiplicities on fibrations  
under birational transformations. The   idea is to use the monoid 
$\SL_2(\NN)$ in its various disguises. In the next
section, they shall apply the results for the classification of singular del Pezzo surfaces.

Rather than working over a ground  field, we shall use
the following set-up: Suppose 
$A$ is a  henselian discrete valuation ring with residue field 
$k$ and fraction field 
$K$, and let 
$X$ be a proper normal 2-dimensional 
$A$-scheme with 
$\Gamma(X,\O_{X})=A$. We call such schemes 
\emph{$A$-surfaces}.

The generic fiber 
$X\otimes K$ is a normal curve   over  the function field
$K$, degenerating into the possibly singular closed fiber 
$X\otimes k$. Call 
$X$ is \emph{minimal} if the closed fiber is irreducible.
If 
$X$ is not minimal, a result of Bosch, L\"utkebohmert and Raynaud 
(\cite{Bosch; Luetkebohmert; Raynaud 1990} Cor.~3, p.~169)
ensures that each irreducible component of   
$X\otimes k$ is contractible (here we need the assumption that 
$A$ is henselian).

Suppose 
$Y$ is  a minimal 
$A$-surface as above. We seek to describe the   passage to birational minimal
$A$-surface algorithmically. Let 
$y\in Y$ be  a closed point in 
$\Reg(Y)$ and 
$f_1\co X_1\ra Y$ be its blowing-up. In the following, 
we consider sequences of blowing-ups 
$$
X_{n+1}\stackrel{f_{n+1}}{\lra} X_{n} \lra\ldots\lra  
X_2\stackrel{f_2}{\lra} X_{1} \stackrel{f_1}{\lra} Y,
$$
subject to the condition that  each center 
$x_i\in X_i$ lies on both the exceptional divisor 
$E_i=f^{-1}_i(x_{i-1})$ of the preceding blowing-up and another 
irreducible component of 
$X_i\otimes k$. This implies that the
intersection graph of 
$X_i\otimes k$ is a chain 
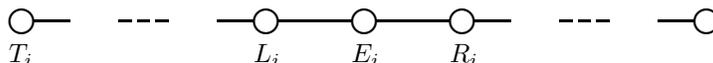
\begin{center}
\setlength{\unitlength}{0.00043300in}%
\begingroup\makeatletter\ifx\SetFigFont\undefined%
\gdef\SetFigFont#1#2#3#4#5{%
  \reset@font\fontsize{#1}{#2pt}%
  \fontfamily{#3}\fontseries{#4}\fontshape{#5}%
  \selectfont}%
\fi\endgroup%
\begin{picture}(8716,900)(1043,-913)
\thicklines
\put(4201,-361){\circle{300}}
\put(5401,-361){\circle{300}}
\put(6601,-361){\circle{300}}
\put(9601,-361){\circle{300}}
\put(1351,-361){\line( 1, 0){450}}
\put(4351,-361){\line( 1, 0){900}}
\put(5551,-361){\line( 1, 0){900}}
\put(4051,-361){\line(-1, 0){450}}
\put(6751,-361){\line( 1, 0){450}}
\put(9451,-361){\line(-1, 0){450}}
\put(1201,-361){\circle{300}}
\put(8401,-361){\line(-1, 0){150}}
\put(6450,-920){\makebox(0,0)[lb]{$R_i$}}
\put(8176,-361){\line(-1, 0){150}}
\put(7951,-361){\line(-1, 0){ 75}}
\put(7951,-361){\line(-1, 0){150}}
\put(2401,-361){\line( 1, 0){150}}
\put(2626,-361){\line( 1, 0){150}}
\put(2851,-361){\line( 1, 0){150}}
\put(1050,-920){\makebox(0,0)[lb]{$T_i$}}
\put(4050,-920){\makebox(0,0)[lb]{$L_i$}}
\put(5250,-920){\makebox(0,0)[lb]{$E_i$}}
\end{picture}
\end{center}
Here 
$T_i$ is the strict transform of 
$Y\otimes k$, and 
$E_i$ is the exceptional divisor of the blowing-up
$f_i\co X_i\ra X_{i-1}$. The two neighbors of 
$E_i$ are labeled
$L_i$ and 
$R_i$; they are distinguished by the fact that 
$L_i$ lies between 
$S_i$ and 
$E_i$.
In the special case 
$i=1$ we have 
$T_1=L_1$ and 
$R_1=\emptyset$. 

Observe that there are two alternatives for the blowing-up 
$f_{i+1}\co X_{i+1}\ra X_i$: Either the center $x_i\in X_i$ is 
$L_i\cap E_i$ or  
$E_i\cap R_i$. In other words: The sequence 
$X_n,\ldots,X_1$ is uniquely determined by the initial center 
$y\in Y$, together with a word 
$w=L^{l_1}R^{r_1}L^{l_2}\ldots$ of length 
$n\geq 0$ in  two letters 
$L,R$ not starting with 
$R$. The rule is: If the 
$i$th letter in 
$w$ is 
$L$,   the center $x_i\in X_i$ is 
$L_i\cap E_i$, otherwise 
$E_i\cap R_i$.

Set $X=X_{n+1}$, let 
$h\co X\ra Y$ be the composition of the blowing-ups,  and 
$\hat{h}\co X\ra \widehat{Y}$ be the contraction of all irreducible components in
$X_{n+1}\otimes k$ except for
$E_{n+1}$. We call the resulting minimal 
$A$-surface $\widehat{Y}$ the \emph{elementary transformation} of 
$Y$ with respect to the word 
$w=L^{l_1}R^{r_1}\ldots$ and the initial center $y\in Y$.

Starting with a regular minimal
$A$-surface, it is easy to see that each birational equivalent minimal  
$A$-surface can be reached by a sequence of elementary transformations.
If 
$Y$ is a 
$\PP^1$-bundle, and
$y$ is a $k$-rational point,  and 
$w=\emptyset$ is the empty word, then the construction
coincides with the classical notion of elementary transformation.

Let   
$\langle L,R\rangle$ be  the free monoid on two letters and   
$\SL_2(\NN)$ the monoid of unimodular
$2\times 2$-matrices whose entries are natural numbers. By inspection, the map 
$$
\langle L,R\rangle \lra \SL_2(\NN),\quad
L\mapsto \begin{pmatrix} 1&1 \\ 0&1 \end{pmatrix},\quad
R\mapsto
\begin{pmatrix} 1&0 \\ 1&1 \end{pmatrix}
$$
is bijective. Similarly,  the map 
$$
\SL_2(\NN)\lra \QQ_+,\quad 
\begin{pmatrix} a&b \\c &d \end{pmatrix}\mapsto \frac{a+b}{c+d}
$$
is bijective. Composing both maps, we associate to each word 
$w=L^{l_1}R^{r_1}\ldots$ a fraction 
$r/s>0$.

Now let 
$\widehat{Y}$ be the elementary transformation of 
$Y$ with respect to the word $w=L^{l_1}R^{r_1}\ldots$
and the initial center $y\in Y$. The   multiplicities 
$\mu,\hat{\mu}>0 $ of the closed fibers are defined by 
$Y\otimes k= \mu \cdot (Y\otimes k)^\red$ and 
$\widehat{Y}\otimes k= \hat{\mu}\cdot (\widehat{Y}\otimes k)^\red$.

\begin{proposition}
\label{multiplicity of fibre}
Let 
$r/s$ be the fraction associated to the word
$w$. Then the multiplicities $\mu,\hat{\mu}$ of the closed fibers of 
$Y,\widehat{Y}$ are related by 
$\hat{\mu}=r\cdot \mu$.
\end{proposition}

\proof
The trick is  to  introduce alternative names for the curves on the blowing-ups. Let 
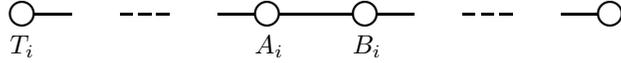
\begin{center}
\setlength{\unitlength}{0.00043300in}%
\begingroup\makeatletter\ifx\SetFigFont\undefined%
\gdef\SetFigFont#1#2#3#4#5{%
  \reset@font\fontsize{#1}{#2pt}%
  \fontfamily{#3}\fontseries{#4}\fontshape{#5}%
  \selectfont}%
\fi\endgroup%
\begin{picture}(7516,1000)(1043,-913)
\thicklines
\put(4201,-361){\circle{300}}
\put(5401,-361){\circle{300}}
\put(8401,-361){\circle{300}}
\put(1351,-361){\line( 1, 0){450}}
\put(4351,-361){\line( 1, 0){900}}
\put(4051,-361){\line(-1, 0){450}}
\put(2401,-361){\line( 1, 0){150}}
\put(2626,-361){\line( 1, 0){150}}
\put(1201,-361){\circle{300}}
\put(2851,-361){\line( 1, 0){150}}
\put(5250,-886){\makebox(0,0)[lb]{$B_i$}}
\put(5551,-361){\line( 1, 0){450}}
\put(8251,-361){\line(-1, 0){450}}
\put(6601,-361){\line( 1, 0){150}}
\put(6826,-361){\line( 1, 0){150}}
\put(7051,-361){\line( 1, 0){150}}
\put(1050,-886){\makebox(0,0)[lb]{$T_i$}}
\put(4050,-886){\makebox(0,0)[lb]{$A_i$}}
\end{picture}
\end{center}
be the intersection graph of 
$X_i\otimes k$. Here 
$T_i$ is the strict transform of 
$Y\otimes k$ as above, but 
$A_i,B_i$ are the curves with 
$A_i\cap B_i =\left\{ x_i \right\}$. They can be distinguished by the fact that 
$A_i$ lies between 
$T_i$ and 
$B_i$.
Let 
$a_i,b_i$ be the multiplicities of 
$A_i,B_i$ in the cycle $X_i\otimes k$. Then 
$\hat{\mu}=a_n+b_n$. 

We proceed to calculate the pairs 
$(a_i,b_i)$. Let 
$a_0=\mu$ and 
$b_0=0$ be dummy variables.
Suppose 
$f_{i+1}\co X_{i+1}\ra X_i$ is defined by the letter 
$L$. Then 
$A_i=L_i$, 
$B_i=E_i$, and 
$$
(a_i,b_i)= (a_{i-1},a_{i-1}+b_{i-1}) = (a_{i-1},b_{i-1})
\begin{pmatrix} 1&1 \\ 0&1 \end{pmatrix}.
$$
Otherwise 
$A_i=E_i$, 
$B_i=R_i$, and 
$$
(a_i,b_i)=(a_{i-1}+b_{i-1},b_{i-1}) = (a_{i-1},b_{i-1})
\begin{pmatrix} 1&0 \\ 1&1 \end{pmatrix}.
$$
Writing
$$
\begin{pmatrix} a&b \\ c&d \end{pmatrix} =
\begin{pmatrix} 1&1 \\ 0&1 \end{pmatrix}^{l_1}
\begin{pmatrix} 1&0 \\ 1&1 \end{pmatrix}^{r_1}
\begin{pmatrix} 1&1 \\ 0&1 \end{pmatrix}^{l_2}\cdots,
$$
we have $r/s=(a+b)/(c+d)$. Inductively, 
$$
(a_n, b_n) = 
(a_0,b_0)\begin{pmatrix} a&b \\ c&d \end{pmatrix}  =
(\mu,0)\begin{pmatrix} a&b \\ c&d \end{pmatrix} =\mu(a, b),
$$
hence 
$\hat{\mu }= a_n+b_n = r\cdot \mu $.
\qed

\medskip
The next task is to calculate the contribution of the relative canonical divisor 
$K_{X/Y}$ on the elementary transformation 
$\widehat{Y}$.

\begin{proposition}
\label{multiplicity of class}
Let 
$r/s$ be the fraction associated to the word
$w$, and 
$\mu$ be the multiplicity of the closed fiber
$Y\otimes k$. Then 
$$
\hat{h}_*(K_{X/Y})= (r+s-1)\cdot( \widehat{Y}\otimes k) ^\red = 
\frac{r+s-1}{r\mu}\cdot (\widehat{ Y}\otimes k).
$$
\end{proposition}

\proof
We use the same notation as in the preceding proof, except that 
$a_i,b_i$ denote the multiplicities of 
$A_i,B_i$ in 
$K_{X_i/Y}$. We introduce dummy multiplicities 
$a_0=0$ and $ b_0=0$. By the adjunction formula, the multiplicity of $E_{n+1}$
in $K_{X_{n+1}/Y}$ is 
$a_n+b_n+1$. It turns out that  the pairs 
$(a_i+1,b_i+1)$ are easy to compute. 

Suppose 
$f_{i+1}\co X_{i+1}\ra X_i$ is defined by the letter 
$L$. Then 
$A_i=L_i$, 
$B_i=E_i$, and 
$$
(a_i+1,b_i+1)= (a_{i-1}+1,a_{i-1}+1+b_{i-1}+1) = (a_{i-1}+1,b_{i-1}+1)
\begin{pmatrix} 1&1 \\ 0&1 \end{pmatrix}.
$$
Otherwise 
$A_i=E_i$, 
$B_i=R_i$, and 
$$
(a_i+1,b_i+1)=(a_{i-1}+1+b_{i-1}+1,b_{i-1}+1) = (a_{i-1}+1,b_{i-1}+1)
\begin{pmatrix} 1&0 \\ 1&1 \end{pmatrix}.
$$
Writing 
$$
\begin{pmatrix} a&b \\ c&d \end{pmatrix} =
\begin{pmatrix} 1&1 \\ 0&1 \end{pmatrix}^{l_1}
\begin{pmatrix} 1&0 \\ 1&1 \end{pmatrix}^{r_1}
\begin{pmatrix} 1&1 \\ 0&1 \end{pmatrix}^{l_2}\cdots,
$$
we have 
$r/s=(a+b)/(c+d)$. Inductively, 
$$
(a_n+1, b_n+1) = 
(a_0+1,b_0+1)\begin{pmatrix} a&b \\ c&d \end{pmatrix}  =
(1,1)\begin{pmatrix} a&b \\ c&d \end{pmatrix} =(a+c, b+d).
$$
Hence the multiplicity of 
$\hat{h}_*(K_{X/Y})$   is
$$
a_n+b_n+1= a+c+b+d-1 = r+s-1.
$$
Together with Proposition \ref{multiplicity of fibre}, 
this concludes the proof.
\qed

\medskip
For your convenience, I recall how  self intersection numbers are related to
continued fractions.
Let
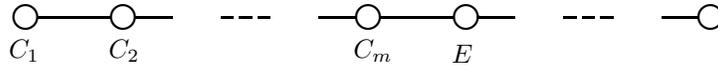
\begin{center}
\setlength{\unitlength}{0.00043300in}%
\begingroup\makeatletter\ifx\SetFigFont\undefined%
\gdef\SetFigFont#1#2#3#4#5{%
  \reset@font\fontsize{#1}{#2pt}%
  \fontfamily{#3}\fontseries{#4}\fontshape{#5}%
  \selectfont}%
\fi\endgroup%
\begin{picture}(8716,934)(1043,-2038)
\thicklines
\put(2401,-1561){\circle{300}}
\put(5401,-1561){\circle{300}}
\put(6601,-1561){\circle{300}}
\put(9601,-1561){\circle{300}}
\put(1351,-1561){\line( 1, 0){900}}
\put(2551,-1561){\line( 1, 0){450}}
\put(3601,-1561){\line( 1, 0){150}}
\put(3826,-1561){\line( 1, 0){150}}
\put(4051,-1561){\line( 1, 0){150}}
\put(5251,-1561){\line(-1, 0){450}}
\put(1201,-1561){\circle{300}}
\put(5551,-1561){\line( 1, 0){900}}
\put(6426,-2111){\makebox(0,0)[lb]{$E$}}
\put(6751,-1561){\makebox(6.6667,10.0000){\SetFigFont{10}{12}{\rmdefault}{\mddefault}{\updefault}.}}
\put(6751,-1561){\line( 1, 0){450}}
\put(7801,-1561){\line( 1, 0){150}}
\put(8026,-1561){\line( 1, 0){150}}
\put(8251,-1561){\line( 1, 0){150}}
\put(9001,-1561){\line( 1, 0){450}}
\put(1000,-2111){\makebox(0,0)[lb]{$C_1$}}
\put(2226,-2111){\makebox(0,0)[lb]{$C_2$}}
\put(5226,-2111){\makebox(0,0)[lb]{$C_m$}}
\end{picture}
\end{center}
be the intersection graph of 
$X\otimes k$, labeled in such a way that 
$C_1$ is the reduced strict transform of 
$Y\otimes k$, and 
$E$ is the reduced strict transform of 
$\widehat{Y}\otimes k$. 

Let 
$r/s>0$ be the fraction associated to the word 
$w=L^{l_1}R^{r_1}\ldots$ describing the elementary transformation, and 
$$
r/s = [s_1,s_2,\ldots, s_m] =  s_1-  \frac{1}{\displaystyle s_2 -\frac{1}{\displaystyle s_3 -\frac{1}{\ldots}}}
$$
its unique continued fraction development with 
$m\geq1$, 
$s_1>0$, and 
$s_i>1$ for 
$i\geq 2$.

\begin{proposition}
\label{selfintersection numbers}
Let 
$r/s = [s_1,s_2,\ldots, s_m]$ be the continued fraction as above. Then 
$C_i^2 = -s_i$. Here the intersection numbers are   computed as degrees over 
$\k(y)$.
\end{proposition}

\proof
The set 
$\foS$ of all finite sequences 
$(s_1,s_2,\ldots, s_m)$  with 
$m\geq1$, 
$s_1>0$, and 
$s_i>1$ for 
$i\geq 2$ becomes a monoid with respect to the composition law
$$
(s_1,\ldots,s_{m-1},  s_m) \circ (t_1,t_2,\ldots, t_n) = (s_1,\ldots, s_{m-1},s_m-1+t_1, t_2,\ldots, t_n).
$$
By inspection, the map
$\langle L,R\rangle \ra \foS$ defined by
$L\mapsto (2)$ and $ R\mapsto (1,2)$ is bijective. We also have a bijective  map 
$\foS\ra\QQ_+$ by 
$(s_1,\ldots, s_m)\mapsto [s_1,\ldots, s_m]$. The maps obtained so far fit into a commutative
diagram 
$$
\begin{CD}
\langle L,R\rangle    @>>>   \SL_2(\NN)  \\
@VVV     @VVV\\
\foS     @>>>   \QQ_+.
\end{CD}
$$
The rest of the proof is much alike the preceding proofs and left as an exercise.
\qed

\medskip
 Let 
$S\subset Y$ be a section disjoint from the initial center 
$y\in Y$, and  
$R\subset X$, $\widehat{S}\subset \widehat{Y}$ its strict transforms. With the notation introduced before the last Proposition,
Mumford's rational pullback takes the form 
$$
\hat{h}^*(\widehat{S}) = R + \gamma_1C_1 +\gamma_2 C_2 +\ldots+\gamma_m C_m
$$
for certain rational coefficients 
$\gamma_i\geq 0$. For future reference, we    determine the first coefficient: 

\begin{proposition}
\label{pullback}
Let 
$r/s$ be the fraction associated to the word 
$w$, and
$\mu$ the multiplicity of 
$Y\otimes k$, and 
$l=\dim_k\k(y)$. Then 
$\gamma_1= s/(r  l\mu)$.
\end{proposition}

\proof
By definition of Mumford's pullback, the coefficients 
$\gamma_i$ solve the equations 
$$
(R + \gamma_1C_1 +\gamma_2 C_2 +\ldots+\gamma_m C_m)\cdot C_i =0,\quad
i=1,\ldots, m.
$$
Observe that $C_1\cdot R= 1/\mu$ and
$R\cap C_i=\emptyset$ for 
$i>1$, since 
$S\subset Y$ is  a section disjoint from the initial center
$y\in Y$. Using Kronecker's 
$\delta$-function, we can rewrite   the equations as
$(\gamma_1C_1   +\ldots+\gamma_m C_m)\cdot C_i =-\delta_{i,1}/\mu$. 
 Introducing dummy variables 
$\gamma_0=1/(l\mu)$ and 
$\gamma_{m+1}=0$, we    put   the equations   into   simplified form
$$
\gamma_{i-1} -s_i\gamma_i +\gamma_{i+1}=0,\quad
i=1,\ldots, m.
$$ 
Here the correction factor 
$l=\dim_k \k(y)$ is necessary since   we calculate   intersection numbers over 
$\kappa(y)$ and not over $k$. We have 
$ \gamma_0/\gamma_1=1/(l\mu\gamma_1)$.
On the other hand, an easy inductive argument gives 
$$
\gamma_0/\gamma_1 = [s_1,\ldots, s_m] = r/s,
$$
and the result follows.
\qed

\section{Classification of del Pezzo surfaces}\label{Classification of del Pezzo surfaces}

In this section, 
$Z$ is a proper normal del Pezzo surface with base number 
$\rho(Z)=1$ containing a nonrational singularity. 
The issue is to understand the geometry
of such surfaces. We   show that these surfaces are arranged in 
an infinite hierarchy;
moving inside the hierarchy involves elementary transformations.
Let 
$f\co X\ra Z$ be the minimal resolution of singularities and 
$$
X\stackrel{h}{\lra} Y \stackrel{g}{\lra}Z
$$
the factorization over the the minimal resolution of nonrational singularities.
According to Theorem \ref{structure}, there is  a fibration 
$p\co Y\ra B$ of genus zero over a curve of genus 
$g>0$ with irreducible fibers, and the exceptional curve 
$S\subset Y$ is a section. 

In the simplest case, 
$Y$ is already a
$\PP^1$-bundle. Otherwise, we seek to relate 
$Y$ with a $\PP^1$-bundle. Let 
$R'\subset X$ be the strict transform of 
$S\subset Y$, and 
$X\ra Y'$ be the contraction of all vertical curves disjoint to 
$R'$. Then 
$Y'$ is smooth, and the induced fibration 
$p'\co Y'\ra B$ is a 
$\PP^1$-bundle. Hence  the normal algebraic  surface
$Y$ is obtained from the smooth surface
$Y'$ by a sequence 
$$
\renewcommand\labelstyle{\scriptstyle}
\begin{diagram}
&X_1&&X_2& &X_{n-2} &&X_{n-1}&\\
\ldTo(1,2)_{h_1}&&\rdTo(1,2)^{\hat{h}_{1}}\ldTo(1,2)_{h_2}&&\cdots&&
\rdTo(1,2)^{\hat{h}_{n-2}}\ldTo(1,2)_{h_{n-1}}&&\rdTo(1,2)^{\hat{h}_{n-1}}\\
\makebox[0.5cm][r]{$Y'=Y_1$}&&Y_2&&&&Y_{n-1}&&\makebox[0.5cm][l]{$Y_n=Y$}\\
\end{diagram}
$$
of elementary transformations. Note that each 
$h_i\co X_i\ra Y_i$ decomposes into a sequence of blowing-ups 
$$
X_i = X_{i, n_i+1}\lra X_{i,n_i} \lra\ldots\lra  
X_{i,2}\lra X_{i,1} \lra Y_i
$$
as in Section \ref{Elementary transformations}.
We have to  fix more notation. Let $R_i\subset X_i$ and 
$S_i\subset Y_i$ be the strict transforms of 
$S\subset Y$.
Let 
$y_i\in Y_i$ be the center of the elementary transformation 
$Y_i\leftarrow X_i\ra Y_{i+1}$,  and   
$w_i\in \langle L,R\rangle$ the corresponding word. Set 
$b_i=p_i(y_i)$, let 
$d_i$ be the 
$k$-dimension of 
$\k(b_i)$, and let 
$r_i/s_i$ be the fraction associated to the word 
$w_i$. Denote by 
$\mu_i$ the multiplicity of 
$Y_i\otimes\k(b_i)$, and by 
$\lambda_i$ the coefficient in 
$K_{Y_i}\equiv -2S_i +\lambda_i F_i$. Here 
$F_i\in N(Y_i)$ is the fiber class. These numbers obey the following rule:

\begin{lemma}
\label{coefficients}
$\lambda_{i+1}= \lambda_i + d_i(r_i+s_i-1)/(\mu_i r_i)$.
\end{lemma}

\proof
By Proposition \ref{multiplicity of fibre}, the multiplicity of the fiber
$Y_{i+1}\otimes \k(b_i)$ is the numerator
$\mu_i r_i$. Moreover, we have 
$Y_{i+1}\otimes \k(b_i)\equiv d_iF_{i+1}$. According to Proposition \ref{multiplicity of class},  the denominator 
$r_i+s_i-1$ is the multiplicity in 
$(\hat{h}_i)_*(K_{X_i/Y_i})$. Since 
$$
 -2S_{i+1} + \lambda_{i +1}F_{i+1} \equiv K_{Y_{i+1}} \equiv -2S_{i+1} + \lambda_i F_{i+1} +
(\hat{h}_i)_*(K_{X_i/Y_i}),
$$
the claim is true.
\qed

\begin{theorem}
\label{structure  of del Pezzo}
For each surface 
$Y_i$, the curve 
$S_i\subset Y_i$ is contractible; the resulting contraction 
$g_i\co Y_i\ra Z_i$ yields a normal del Pezzo surface 
$Z_i$ with base number
$\rho(Z_i)=1$ containing a nonrational singularity.
\end{theorem}

\proof
We verify the assertion by descending induction on 
$i$. For 
$i=n$, there is nothing to prove. Suppose that the statement is true for some 
$i+1\leq n$. It follows from the definition of Mumford's rational intersection numbers that 
$R_i^2\leq S_{i+1}^2<0$. Since 
$X_i\ra Y_i$ is locally an isomorphism near 
$S_i$, we have 
$S_i^2=R_i^2$. Thus 
$S_i$ has negative selfintersection.
By induction, we have 
$\lambda_{i+1}<0$. According to Lemma \ref{coefficients},
$\lambda_i\leq \lambda_{i+1}$. 
Hence  Proposition \ref{contractible curve} applies, and the contraction 
$g_i\co Y_i\ra Z_i$ of 
$R_i\subset Y_i$ exists. Since the canonical class is given by
$K_{Z_i}= \lambda_i\cdot (g_i)_*(F_i)$, the resulting surface is del Pezzo.
\qed

\medskip
We see that  our del Pezzo surface 
$Z=Z_n$ is the last in a sequence 
of del Pezzo surface $Z_1,\ldots, Z_n$  underlying a sequence of    
elementary transformations
$$
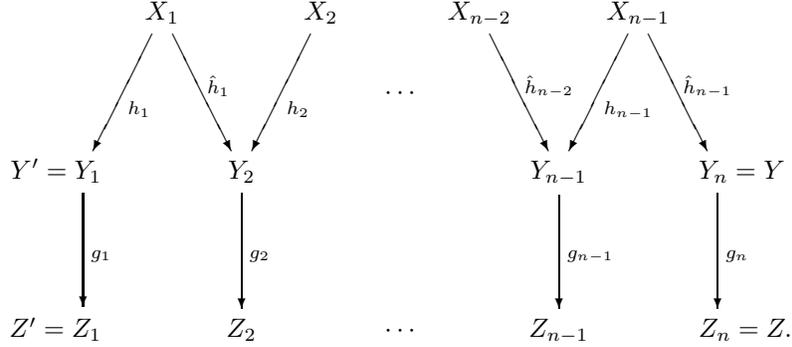
\renewcommand\labelstyle{\scriptstyle}
\begin{diagram}
&X_1&&X_2& &X_{n-2} &&X_{n-1}&\\
\ldTo(1,2)_{h_1}&&\rdTo(1,2)^{\hat{h}_{1}}\ldTo(1,2)_{h_2}&&\cdots&&
\rdTo(1,2)^{\hat{h}_{n-2}}\ldTo(1,2)_{h_{n-1}}&&\rdTo(1,2)^{\hat{h}_{n-1}}\\
\makebox[0.5cm][r]{$Y'=Y_1$}&&Y_2&&&&Y_{n-1}&&\makebox[0.5cm][l]{$Y_n=Y$}\\
\dTo>{g_1}&&\dTo>{g_2}&&&&\dTo>{g_{n-1}}&&\dTo>{g_n}\\
\makebox[0.5cm][r]{$Z'=Z_1$}&&Z_2&&\cdots&&Z_{n-1}&&\makebox[0.5cm][l]{$Z_n=Z.$}\\
\end{diagram}
$$
The initial term 
$Z'=Z_1$ is a 
contraction of the 
$\PP^1$-bundle $Y'=Y_1$.
To obtain  complete command of the situation, 
we have to describe how such   sequence  comes about.
The question is: can we prolong the sequence one step further?

Suppose
$Y_n\leftarrow X_n\ra Y_{n+1}$ is a   elementary transformation with center 
$y_n\in Y_n$ in 
$\Reg(Y_n)$ disjoint from $S_n\subset Y_n$.
Introduce the same notation laid down in front of 
Theorem \ref{structure of del Pezzo} as for the preceding elementary
transformations. The issue is to determine whether or not the strict transform 
$S_{n+1}\subset Y_{n+1}$  contracts to a del Pezzo surface.

\begin{theorem}
\label{existence of del Pezzo}
Notation as above. Then the curve $S_{n+1}\subset Y_{n+1}$  
contracts to a del Pezzo surface
$Z_{n+1}$ if and only if  
$$
\lambda_n + d_n\frac{r_n+s_n-1}{\mu_n r_n} <0.
$$
\end{theorem}

\proof
According to Lemma \ref{coefficients}, 
$\lambda_{n+1}= \lambda_n + d_n(r_n+s_n-1)/(\mu_n r_n)$.
This being the coefficient in 
$K_{Y_{n+1}}\equiv -2S_{n+1} +\lambda_{n+1} F_{n+1}$,  
the condition  is necessary.

Conversely, assume that the condition holds. Then  
$\lambda_{n+1}< 0$.
The main problem is to check that 
$S_{n+1}\subset Y_{n+1}$ has negative selfintersection.
I claim that 
$S_i^2\leq \lambda_i$ holds for 
$1\leq i\leq n+1$. Suppose for a moment that this is true. Then 
$S_{n+1}$ has negative selfintersection, and the desired contraction exists 
by to Proposition
\ref{contractible curve}.

It remains to prove the claim. We proceed by induction on 
$i$. For 
$i=1$, the fibration 
$p_1\co Y_1\ra B$ is a 
$\PP^1$-bundle. The adjunction formula gives 
$$
K_{Y_1} \equiv -2S_1 + (S_1^2 + 2g-2)F_1,
$$
hence 
$S_1^2\leq S_1^2 + 2g-2 = \lambda_1$.

Suppose the claim is correct for some 
$i\geq 1$. Let 
$T_i\subset X_i$ be the reduced strict transform of the fiber
$Y_{i+1}\otimes \k(b_{i+1})$. Write 
$$
(\hat{h}_{i})^*(S_{i+1}) = R_i + \gamma_i \cdot T_i +\ldots$$
for some rational coefficient 
$\lambda_i\geq 0$. With 
$l_i=\dim_{\k(b_i)} \k(y_i)$,  Proposition \ref{pullback} gives 
$\gamma_i = s_i/(r_i l_i\mu_i)$.
Note that 
$R_i\subset X_i$ hits the fiber over 
$b_i$ only in 
$T_i$, since the center 
$y_i\in Y_i$ is disjoint from 
$S_i$. Using Mumford's definition of rational pullback, we obtain 
\begin{eqnarray*}
S_{i+1}^2 & =&  R_i\cdot (R_i + \gamma_i T_i +\ldots) \\
                  & = & R_i^2 + d_i/\mu_i\cdot s_i/ (r_i l_i\mu_i)\\
                  & = & S_i^2 + d_is_i/ (r_i l_i\mu_i^2)\\
                  &\leq& \lambda_i + d_i(s_i+r_i-1)/ (r_i\mu_i).
\end{eqnarray*}
By Lemma \ref{coefficients}, the latter expression equals
$\lambda_{i+1}$. This completes the induction, and therefore the whole proof.
\qed

\section{Anticanonical rings and anticanonical models}
\label{anticanonical model}

Each proper normal algebraic surface 
$X$ comes along with the \emph{anticanonical ring} 
$$
R(-K_X)=\bigoplus_{n\geq 0}H^0(X, \O_{X}(-nK_X)),$$
which in turn defines the 
\emph{anticanonical model} 
$P(-K_X) =\Proj R(-K_X)$.
In this section we shall study the scheme 
$P(-K_X)$.
Zariski (\cite{Zariski 1962}  p.~562) observed that the ring
$R(-K_X)$ is possible non-Noetherian. Surprisingly, the scheme 
$P(-K_X)$ is  of finite type (\cite{Schroeer 1999a} Thm.\ 6.2). 
It is either empty, a point, a proper normal curve,
or a normal algebraic  surface (possibly nonproper). 

There is   a rational map 
$r\co X\dasharrow P(-K_X)$ defined as follows: 
By definition, the homogeneous spectrum 
$P(-K_X)$ is covered by affine open subsets 
$D_+(s)=\Spec R_{(s)}$ with 
$s\in H^0(-nK_X)$ and
$n>0$. Let 
$X_s\subset X$ be the open subset where 
$s\co \O_{X}\ra\O_{X}(-nK_X)$ is bijective. Then the affine hull 
$X_s^\aff=\Spec\Gamma(X_s,\O_{X})$ and 
$D_+(s)$ are canonically isomorphic, and we obtain a morphism 
$r\co U\ra P(-K_X)$ defined on 
$U=\bigcup X_s$. We call 
$\SBs(-K_X) = X\setminus U$ the  \emph{stable base locus}.

\begin{theorem}
\label{anticanonical}
Let 
$X$ be a proper normal algebraic  surface with 2-dimensional anticanonical model 
$P=P(-K_X)$. Then there is a unique   open dense embedding 
$P\subset \P$  into a proper normal  algebraic surface 
$\P$  such that the following holds:
\renewcommand{\labelenumi}{(\roman{enumi})}
\begin{enumerate}
\item The open subset
$P\subset \P$ is the 
$\QQ$-Gorenstein locus of $\P$.
\item The boundary at infinity 
$\P\setminus P$ contains at most one point.
\end{enumerate}
\end{theorem}

\proof
Choose some integer
$n>0$ such that the stable base locus 
$\SBs(-K_X)$ is the base locus of 
$-nK_X$. Let 
$F\subset X$ be the fixed curve of 
$-nK_X$, and let
$M=-nK_X-F$ be the movable part.
Since 
$P$ is 2-dimensional, we may furthermore assume that 
$M\neq 0$, such that 
$h^0(-nK_X)>1$.
Let $F'\subset F$   be the union of all connected components 
$C\subset F$ with 
$M\cdot C>0$, and set 
$F''=F-F'$.

The idea is to construct the compactification $\P$ as a suitable
contraction of $X$.
To do so, we have to define the corresponding negative definite curves.
Let 
$R\subset X$ be the union of all curves  
$C\subset X$ with 
$C\cap F'=\emptyset$ and 
$C\cdot M=0$. Note that 
$F''\subset X$. By the Hodge index theorem, the curve
$R$ is negative definite. 

\begin{claim}
The curve $F'\subset X$ is also negative definite.
\end{claim}

To prove this, let 
$C\subset F'$ be  a connected component  and 
$F_1+\ldots+F_r$ be its integral components. We have to show that the intersection matrix 
$(F_i\cdot F_j)$ is negative definite. Seeking a contradiction, we first assume that some Weil divisor 
 supported by 
$C$  has positive selfintersection.  As in the proof of \cite{Schroeer 1999b} Proposition 3.2, we find a
curve 
$A\subset X$ with 
$\Supp(A)=\Supp(C)$, such that 
$A\subset X$ is linear equivalent to a curve 
$B\subset X$ not containing any 
$F_i$. Decompose 
$A=\sum\lambda_iF_i$ and 
$C=\sum\mu_iF_i$. We may assume that 
$\lambda_1/\mu_1\geq\lambda_i/\mu_i$ for all indices
$1\leq i\leq r$. Then 
$$
\lambda_1 C = \sum_{i=1}^{r}\lambda_1\mu_iF_i =  
\sum_{i=1}^{r}\mu_1\lambda_iF_i
+ \sum_{i=2}^{r}(\lambda_1\mu_i-\lambda_i\mu_1)F_i,
$$
which is linearly equivalent to the effective divisor
$\mu_1B+\sum_{i=2}^{r}(\lambda_1\mu_i-\lambda_i\mu_1)F_i$ not containing 
$F_1$. Consequently, 
$F_1\not\subset\SBs(-K_X)$, contradiction.

Next, assume that 
$(F_i\cdot F_j)$ is negative, but not definite. 
By \cite{Barth; Peters; Van de Ven 1984} Lemma I.2.10, the intersection form 
is negative
semidefinite, and we  find a curve 
$A\subset X$ with 
$\Supp(A)=\Supp(C)$ and 
$A\cdot F_i=0$ for  
$1\leq i\leq r$. Riemann--Roch for normal surfaces gives 
$$
h^0(tA)+h^2(tA) \geq 
\frac{A^2}{2}t^2 -\frac{A\cdot K_X}{2} t+\chi(\O_{X}) +\psi(t)
$$
for some bounded error function 
$\psi(t)$.
Since 
$K_X\cdot A<0$ and 
$K_X-tA$ is not effective,  we conclude that some 
$tA$ is linearly equivalent to a curve 
$B\subset X$ disjoint from 
$C$. As above, this contradicts
$C\subset\SBs(-K_X)$. QED for the claim.

\medskip
We proceed with the proof of the Theorem.
By Proposition \ref{contractions on del Pezzo}, the negative definite curve 
$F'\cup R\subset X$ is contractible. Let 
$f\co X\ra Z$ be its contraction. We  show that 
$Z=\P$ is the desired compactification of the anticanonical model 
$P=P(-K_X)$.  We have to construct an open embedding 
$P\subset Z$. Given a section
$s\in H^0(-nK_X)$, let 
$t\in H^0(M+F')$ be the corresponding section. You easily check that 
$D_+(s)=X_s^\aff$ is isomorphic to the open subset 
$X_t^\aff\subset Z$. Patching these isomorphisms gives the desired 
open embedding 
$P\subset Z$.
By construction, the boundary at infinity $Z\setminus P$ is the image of 
$\SBs(-K_X)\setminus F''$, which is also 
$\SBs(-K_Z)$.

Note that 
$-nK_Z=f_*(-nK_X)=f_*(M)$ is movable, so 
$\SBs(-K_Z)$ is discrete. 
First, suppose that 
$Z$ is 
$\QQ$-Gorenstein. By the Fujita--Zariski theorem 
(\cite{Fujita 1983} Thm.\ 1.10),  the stable base locus 
$\SBs(-K_Z)$ is empty, so 
$Z=P$.
Second, suppose that 
$Z$ contains a non-$\QQ$-Gorenstein point 
$z\in Z$. Since rational surface singularities are 
$\QQ$-factorial, the point 
$z\in Z$ is a nonrational singularity. By Theorem \ref{structure}, 
all other singularities are rational.
Applying the Fujita--Zariski Theorem on the resolution of 
$z\in Z$, we deduce that 
$\SBs(-K_X) =\left\{  z \right\}$. Hence we have a disjoint union
$Z=P\cup \left\{  z \right\}$. In either case, the open subset 
$P\subset Z$ is the 
$\QQ$-Gorenstein locus.
\qed

\medskip
Here is a sufficient condition for 2-dimensional anticanonical models:

\begin{proposition}
\label{anticanonical surface}
If 
$K_X^2>0$ and 
$K_X$ is not pseudoeffective, then the anticanonical model
$P(-K_X)$ is a  surface.
\end{proposition}

\proof
Riemann-Roch for normal surfaces \cite{Giraud 1982}  and Serre duality yield 
$$
h^0(-nK_X) + h^0((1+n)K_X) \geq K^2_X\cdot n^2 +\chi(\O_{X}) + \psi(n)
$$
for certain   bounded error function $\psi(n)$. Hence there is an integer 
$n>0$ and a curve 
$C\subset X$ representing
$-nK_X$. By \cite{Schroeer 1999b} Proposition 3.2, the complement 
$U=X-C$ has   2-dimensional global section ring 
$\Gamma(U,\O_{X})$, since 
$C^2>0$. Because the affine hull
$U^\aff=\Spec\Gamma(U,\O_{X}) $ is an open subset of 
$P(-K_X)$, the assertion follows.
\qed

\begin{example}
Here are normal del Pezzo surfaces with non-finitely generated
canonical rings, which also occur in
\cite{Badescu 1983}. As in Section \ref{contractions of ruled},
let $B$ be a curve of genus $g\geq 2$, and $D\in\Div(B)$ a 
divisor of degree $0<d<2g-2$,
and $Y=\PP(\O_B\oplus\O_B(D))$ the corresponding ruled surface.
According to Proposition \ref{del Pezzo contraction}, 
the negative section $S\subset Y$ is contractible,
and the corresponding  contraction $g\co Y\ra X$ gives a del Pezzo surface.
Let $x\in X$ be the resulting singular point.
You easily check that $\O_{X,x}$ is $\QQ$-Gorenstein if and only if 
the divisors
$K_B$ and  $D$ are linearly dependent in $\Pic(B)\otimes\QQ$.
In this  case, the anticanonical model is
$P(-K_X)=X$, and the anticanonical ring is finitely
generated. Otherwise, we have $P(-K_X)=X-\left\{x\right\}$, 
such that the anticanonical ring
is not finitely generated.
\end{example}



\begin{thebibliography}{ccccc}

\bibitem{Artin 1962}
M.~Artin:
Some numerical criteria for contractibility of curves on algebraic
surfaces.
Am.\ J.\ Math.\ 84, 485--496 (1962).

\bibitem{Badescu 1983}
L.~B\u{a}descu:
Anticanonical models of ruled surfaces. 
Ann.\ Univ.\ Ferrara, Nuova Ser.\   29, 165--177 (1983).

\bibitem{Barth; Peters; Van de Ven 1984}
W.~Barth, C.~Peters, A.~Van de Ven:
Compact complex surfaces. 
Ergeb.\ Math.\  Grenzgebiete (3) 4, 
Springer, Berlin, 1984.

\bibitem{Bosch; Luetkebohmert; Raynaud 1990}
S.~Bosch, W.~L\"utkebohmert, M.~Raynaud:
N\'eron models.  
Ergeb.\ Math.\ Grenzgebiete (3) 21.
Springer, Berlin, 1990.

\bibitem{Cheltsov 1997}
I.~Cheltsov:
Del Pezzo surfaces with nonrational singularities.
Math.\ Notes 62, No.3, 377--389 (1997).

\bibitem{Fujisawa 1995}
T.~Fujisawa:
On non-rational numerical Del Pezzo surfaces. 
Osaka J.\ Math.\ 32, 613--636 (1995).

\bibitem{Fujita 1983}
T.~Fujita:
Semipositive line bundles. 
J.\ Fac.\ Sci.\ Univ.\ Tokyo   30, 353--378 (1983).

\bibitem{Giraud 1982}
J.~Giraud:
Improvement of Grauert-Riemenschneider's theorem for a normal surface.
Ann.\ Inst.\ Fourier 32,  13--23 (1982).

\bibitem{SGA 6}
A.~Grothendieck et al.:
Th\'eorie des intersections et th\'eor\`eme de Riemann-Roch.
Lect.\ Notes   Math.\ 225. 
Springer, Berlin, 1971.

\bibitem{Hartshorne 1970}
R.~Hartshorne:
Ample subvarieties of algebraic varieties.
Lect.\ Notes  Math.\ 156.
Springer, Berlin, 1970.

\bibitem{Hidaka; Watanabe 1981}
F.~Hidaka, K.~Watanabe:
Normal Gorenstein surfaces with ample anti-canonical divisor. 
Tokyo J.\ Math.\ 4, 319--330  (1981).

\bibitem{Hidaka; Tomari 1989}
F.~Hidaka, M.~Tomari:
On singularities arising from the contraction of the minimal section of a
ruled surface.
Manuscr.\ Math.\ 65, 329--347 (1989).

\bibitem{Kleiman 1966}
S.~Kleiman:
Toward a numerical theory of ampleness.
Ann.\ Math.\ (2) 84, 293--344  (1966).

\bibitem{Kollar 1995}
J.~Koll\'ar:
Rational curves on algebraic varieties. 
Ergeb.\ Math.\  Grenzgebiete (3)   32. 
Springer, Berlin, 1995.

\bibitem{Mori 1982}
S.~Mori:
Threefolds whose canonical bundles are not numerically effective. 
Ann.\ Math.\ 116, 133--176 (1982).

\bibitem{Mumford 1961}
D.~Mumford:
The topology of an normal surface singularity of an algebraic variety and criterion for simplicity. 
Publ.\ Math., Inst.\ Hautes \'Etud.\ Sci.\ 9, 5--22 (1961).

\bibitem{Mumford 1966}
D.~Mumford:
Lectures on curves on an algebraic surface.
Princeton University Press, Princeton, 1966.

\bibitem{Mumford 1969}
D.~Mumford:
Enriques' classification of surfaces in char p  I. 
In: D.~Spencer, S.~Iyanaga (eds.), Global Analysis, pp.\ 325--339.
Princeton University Press, Princeton, 1969.

\bibitem{Mumford 1970}
D.~Mumford:
Abelian varieties. 
Tata Institute of Fundamental Research Studies in Mathematics 5.
Oxford University Press,  London, 1970.

\bibitem{Sakai 1984}
F.~Sakai:
Anticanonical models of rational surfaces. 
Math.\ Ann.\ 269, 389--410 (1984).

\bibitem{Sakai 1988}
F.~Sakai:
Ruled fibrations on normal surfaces. 
J.\ Math.\ Soc.\ Japan 40,  249--269 (1988).

\bibitem{Schroeer 1999a}
S.~Schr\"oer:
On non-projective normal surfaces.
Manuscr.\ Math.\ 100, 317--321 (1999).

\bibitem{Schroeer 1999b}
S.~Schr\"oer:
On contractible curves on normal surfaces.
J.\ Reine Angew.\ Math.\ 524, 1--15 (2000).

\bibitem{Tango 1972}
H.\ Tango:
On the behavior of extensions of vector bundles under the Frobenius map. 
Nagoya Math.\ J.\ 48,  73--89 (1972).

\bibitem{Tate 1952}
J.~Tate:
Genus change in inseparable extensions of function fields.
Proc.\ Am.\ Math.\ Soc.\ 3, 400--406 (1952).

\bibitem{Zariski 1962}
O.~Zariski:
The theorem of Riemann-Roch for high multiples 
of an effective divisor on an algebraic surface.
Ann.\ Math.\ (2) 76, 560--615 (1962).

\end{thebibliography}
\end{document}